\immediate\write18{makeindex -s nomencl.ist -o "\jobname.nls" "\jobname.nlo"}
\documentclass{IEEEtran4PSCC}
\usepackage{tikz}
\usetikzlibrary{arrows.meta, positioning, backgrounds, automata, decorations, shapes, fit}
\usepackage{pgfplots}
\usepackage[outline]{contour}
\contourlength{1.2pt}

\tikzset{
    >=stealth,
}
\usepgfplotslibrary{fillbetween}
\pgfplotsset{compat = newest}
\usepgfplotslibrary{polar}
\usepgflibrary{shapes.geometric}
\usetikzlibrary{calc}
\pgfplotsset{colormap/violet}

\usepackage{pdfrender}
\usepackage{comment}
\usepackage{multirow}
\usepackage{booktabs}

\usepackage{subcaption}

\usepackage{siunitx}
\usepackage{xcolor}
\usepackage{hyperref}   

\usepackage{mathtools}
\usetikzlibrary{shapes, arrows, snakes, calc}

\usepackage{amsmath}
\usepackage{graphicx}
\usepackage{latexsym}
\usepackage{color}
\usepackage{amssymb}
\usepackage{accents}
\usepackage{tabularx}
\usepackage{fancyhdr}
\usepackage{verbatim}
\usepackage{multirow}
\usepackage{framed}

\usepackage{url}
\usepackage{float}
\usepackage{enumitem}
\usepackage{mathtools}
\usepackage{mathrsfs}
\usepackage{amsfonts}
\usepackage{listings}
\usepackage{amsthm}
\usepackage{grffile}
\usepackage{sidecap}
\usepackage{pbox}
\usepackage[english,ruled,vlined,linesnumbered,algo2e]{algorithm2e}
\usepackage[noend]{algpseudocode}
\usepackage{scalefnt}
\usepackage{nomencl}
\usepackage{etoolbox}

\makenomenclature
\renewcommand\nomgroup[1]{%
    \item[\bfseries
    \ifstrequal{#1}{A}{Indices and index sets}{%
    \ifstrequal{#1}{B}{Parameters}{%
    \ifstrequal{#1}{C}{Decision variables}{%
    }}}%
]\vspace{0.1in}}

\usepackage{diagbox}
\usepackage{cuted}
\usepackage{array}
\usepackage{textcomp}
\usepackage{stfloats}
\usepackage{bm}
\usepackage{algorithm}

\newtheorem{theorem}{Theorem}

\newtheorem{remark}{Remark}

\newcommand{\cG}{{\cal G}}
\newcommand{\cB}{{\cal B}}
\newcommand{\cL}{{\cal L}}
\newcommand{\cD}{{\cal D}}

\newcommand{\cT}{{\cal T}}

\newcommand{\cC}{{\cal C}}

\newcommand{\st}{\mbox{s.t.}}

\allowdisplaybreaks

\renewcommand{\underbar}{\underaccent{\bar}}
\newcommand{\tcr}{\textcolor{black}}

\newfloat{model}{thp}{lop}
\floatname{model}{Model}

\def\BibTeX{{\rm B\kern-.05em{\sc i\kern-.025em b}\kern-.08em
    T\kern-.1667em\lower.7ex\hbox{E}\kern-.125emX}}
\usepackage{balance}
\makeatletter
\let\old@ps@headings\ps@headings
\let\old@ps@IEEEtitlepagestyle\ps@IEEEtitlepagestyle
\def\psccfooter#1{%
    \def\ps@headings{%
        \old@ps@headings%
        \def\@oddfoot{\strut\hfill#1\hfill\strut}%
        \def\@evenfoot{\strut\hfill#1\hfill\strut}%
    }%
    \def\ps@IEEEtitlepagestyle{%
        \old@ps@IEEEtitlepagestyle%
        \def\@oddfoot{\strut\hfill#1\hfill\strut}%
        \def\@evenfoot{\strut\hfill#1\hfill\strut}%
    }%
    \ps@headings%
}
\makeatother

\psccfooter{%
        \parbox{\textwidth}{\hrulefill \\ \small{23rd Power Systems Computation Conference} \hfill \begin{minipage}{0.2\textwidth}\centering \vspace*{4pt} \includegraphics[scale=0.06]{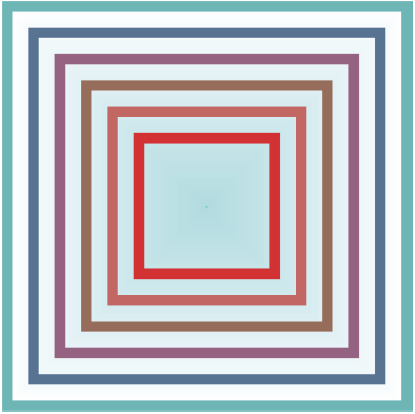}\\\small{PSCC 2024} \end{minipage} \hfill \small{Paris, France --- June 4 -- 7, 2024}}%
}

\begin{document}
%
\title{Multistage Stochastic Program for Mitigating Power System Risks under Wildfire Disruptions}

\author{
\IEEEauthorblockN{Hanbin Yang\IEEEauthorrefmark{1}, Haoxiang Yang\IEEEauthorrefmark{1}}
\IEEEauthorblockA{School of Data Science\\
The Chinese University of Hong Kong, Shenzhen,\\
Shenzhen, China \\
hanbinyang@link.cuhk.edu.cn,\\ yanghaoxiang@cuhk.edu.cn}
\and
\IEEEauthorblockN{Noah Rhodes\IEEEauthorrefmark{2}, Line Roald\IEEEauthorrefmark{2}}
\IEEEauthorblockA{College of Engineering\\
University of Wisconsin, Madison,\\
Madison, WI\\
\{nrhodes, roald\}@wisc.edu}
\and
\IEEEauthorblockN{Lewis Ntaimo\IEEEauthorrefmark{3}}
\IEEEauthorblockA{Industrial \& Systems Engineering\\
Texas A\&M University,\\
College Station, TX\\
ntaimo@tamu.edu}

}

\maketitle

\begin{abstract}
    The frequency of wildfire disasters has surged fivefold in the past 50 years due to climate change. Preemptive de-energization is a potent strategy to mitigate wildfire risks but substantially impacts customers. We propose a multistage stochastic programming model for proactive de-energization planning, aiming to minimize economic loss while accomplishing a fair load delivery. We model wildfire disruptions as stochastic disruptions with varying timing and intensity, introduce a cutting-plane decomposition algorithm, and test our approach on the RTS-GLMC test case. Our model consistently offers a robust and fair de-energization plan that mitigates wildfire damage costs and minimizes load-shedding losses, particularly when pre-disruption restoration is considered.
\end{abstract}


\thanksto{\noindent Submitted to the 23rd Power Systems Computation Conference (PSCC 2024). Hanbin Yang and Haoxiang Yang are supported by the National Natural Science Foundation of China (grant number 72201232) and Shenzhen Key Laboratory of Crowd Intelligence Empowered Low-Carbon Energy Network (project number ZDSYS20220606100601002). Noah Rhodes and Line Roald are supported by the U.S. National Science Foundation under the NSF CAREER award 2045860.}

\section{Introduction}
    \noindent The electric power infrastructure has frequently ignited highly destructive wildfires~\cite{Muhs2020}, leading to catastrophic fires and causing significant loss of life and property. Effective de-energization mitigates wildfire risk of power line ignitions~\cite{rhodes2020balancing}, but excessive de-energization can result in massive load shedding and significant economic and societal impacts~\cite{PacificGas2}. Thus, optimizing de-energization decisions is crucial for balancing wildfire risk mitigation and power outage impacts.

    Previous works primarily address wildfire risks from power system components when optimizing de-energization decisions, overlooking simultaneous natural/human-made wildfire damage to power systems~\cite{rhodes2020balancing} \tcr{and primarily concentrating on the susceptibility of components to wildfires~\cite{taylor2022framework, muhs2020wildfire, davoudi2021reclosing}. Furthermore, such approaches are often deterministic due to computational complexities in solving uncertainty-based optimization models with binary de-energization variables~\cite{vazquez2022wildfire, huang2023review, bayani2022quantifying}. However, incorporating wildfire spread dynamics can significantly enhance de-energization performance. Recent models have integrated stochastic wildfire risks, like dynamic programming for Public Safety Power Shut-offs (PSPS) optimization~\cite{DP_PSPS} and two-stage models for \tcr{mitigating wildfire risk~\cite{yang2023multi, zhou2023optimal, bayani2023resilient}}. Additionally, Ref.~\cite{Noah2022, kody2022sharing} have considered joint de-energization and restoration operations. Despite these developments, existing literature often exhibits a relatively simple decision structure, failing to fully capture the temporal dynamics of wildfires.} Nonetheless, much of the literature employs relatively simple decision structures, failing to capture wildfire's temporal dynamics. Consequently, exploring the interaction between power system operations and dynamic wildfire progress is necessary~\cite{WildfireManagement2}.
        
    Our prior work~\cite{yang2023multi} introduces a two-stage stochastic program but without restoration planning or fairness in power supply allocation. In this paper, we address these limitations by introducing a scenario-based multistage stochastic mixed-integer program (M-SMIP)\tcr{, a framework used to model long-term planning applications involving stochastic uncertainty such as unit commitment~\cite{takriti1996stochastic}, capacity expansion~\cite{rajagopalan1998capacity}, and generation scheduling in hydro systems~\cite{flach2010long}.} We extend the two-stage model to encompass dynamic multistage scenarios and explicitly include restoration considerations within the subproblem to enhance operational resilience. We incorporate fairness constraints in the first problem to ensure equitable power supply allocation to obtain a balanced PSPS plan, as detailed in~\cite{PSPS}. We assess the advantages of incorporating restoration operations while considering fairness, which raises important questions regarding system flexibility and operational insights.

    We propose a decomposition algorithm for the non-convex and non-smooth M-SMIP. Our algorithm is a deterministic variant of SDDP~\cite{pereira1991}, employing linear cutting planes to approximate the value function, thus demonstrating computational precision and reliability.
    
    The contributions of our paper are threefold:
    \begin{enumerate}
        \item We formulate an M-SMIP model that manages wildfire risk in power system operations. Unlike the existing literature, our model captures the random nature of wildfire ignitions using disruption scenarios with random onset times to keep the model size manageable. \tcr{Compared with our prior work in \cite{yang2023multi}, we assume that there multiple disruptions may occur.}
        \item We develop a decomposition algorithm with finite convergence for solving the presented model. We conduct a comparative analysis and investigation of the effectiveness of various cut families, which shows the efficiency of Lagrangian cuts and enhancement.
        \item We show that our M-SMIP model with restoration and fairness yields resilient solutions on the RTS-GLMC test case~\cite{Kody2022}. \tcr{We gain insights from restoration operations, which demonstrate that these operations are superior in reducing wildfire risk and ensuring load delivery.}
    \end{enumerate}
    The remainder of this paper is organized as follows. Section~\ref{sec:model} formulates the optimized PSPS problem. Section~\ref{sec:alg} describes the decomposition algorithm and cut families. Section~\ref{sec:numerical} describes our test case setup and provides numerical results. Section~\ref{sec:conclusion} concludes the paper and discusses future work.

\section{Problem Formulation} \label{sec:model}
    \nomenclature[A,01]{$\cB, \cG, \cL$}{set of buses, generators, and transmission lines;}
    \nomenclature[A,03]{$\cD$}{set of load demand;}
    \nomenclature[A,02]{$\cC$}{set of load components, \(\cC = \cB \cup \cG \cup \cL\);}
    \nomenclature[A,04]{$\cT$}{set of time periods, \(\cT = \{1,2,\dots,T\}\);}
    \nomenclature[A,05]{${\Omega}$}{the set of realizations of wildfire random variables;}

    \nomenclature[B,01]{$D_{dt}$}{load \(d \in \cD\) at time period \(t \in \cT\);}
    \nomenclature[B,02]{$w_{d}$}{the priority level of load \(d \in \cD\);}
    \nomenclature[B,03]{$c^r_{c}$}{repair cost for component \(c \in \cC \);}  
    \nomenclature[B,04]{$\underbar{P}^G_g, \bar{P}^G_g$}{maximum and minimum generation limits of \(g \in \cG \);}  
    \nomenclature[B,05]{$W_{ij}$}{the thermal power flow limit of the line \((i,j) \in \cL \);}  
    \nomenclature[B,06]{$b_{ij}$}{the susceptance of the line \( (i,j) \in \cL \);}  
    \nomenclature[B,07]{$\underbar{\theta}, \bar{\theta}$}{the big-M values for voltage angle difference;} 
    \nomenclature[B,11]{$\beta$}{fairness level;}

    \nomenclature[C,01]{$\theta_{it}^\omega$}{phase angle of the bus \(i \in \cB\) at time \(t \in \cT \) for realization $\omega$;}
    \nomenclature[C,02]{$P^{L,\omega}_{ijt}$}{active power flow on the line \((i,j) \in \cL\) at time \(t \in \cT \) for realization $\omega$;}
    \nomenclature[C,03]{$p^{G,\omega}_{gt}$}{active power generation at generator \(g \in \cG\) at time \(t \in \cT \) for realization $\omega$;}
    \nomenclature[C,04]{$s^\omega_{dt}$}{percentage of load-shedding at the load \(d \in \cD\) at time \(t \in \cT \) for realization $\omega$;}
    \nomenclature[C,06]{$z_{ct}^\omega$}{$1$ if component $c\in \cC$ is functional at time $t \in \cT$ for realization $\omega$, $0$ otherwise;}
    \nomenclature[C,07]{$\nu^\omega_{ct}$}{$1$ if a fire damage is incurred at component $c \in \cC$ at time $t \in \cT$ for realization $\omega$, $0$ otherwise;}
    \nomenclature[C,8]{$r_{ct}$}{$1$ if restoration has not been applied to component $c\in \cC$ for the first-stage problem, $0$ otherwise.}
    \printnomenclature

    \noindent We consider multi-period dispatch and de-energization operations under wildfire disruption over a short-term time horizon (24 hours). The power network is represented by a graph $(\cB, \cL)$, where $\cB$ is the set of buses and $\cL$ is the set of lines. We use $\mathcal{D}_i$, $\mathcal{G}_i$, and $\cL_i$ to represent the subset of loads, generators, and transmission lines connected to bus $i$, respectively. For each period $t \in \mathcal{T}$, we formulate DC power flow constraints to decide the active power generation \(P^G_{gt}\) of generator \(g \in \cG\), the power flow $P^L_{ijt}$ on transmission line \((i,j) \in \cL\) from bus \(i\) to bus \(j\), and the phase angle $\theta_{it}$ of bus \(i \in \cB\).

    \subsection{Wildfire Scenarios}
\begin{figure}
            \centering
            \begin{tikzpicture}[scale = .7, square/.style={regular polygon,regular polygon sides=4},node distance={20mm}, thick, main/.style = {draw, circle}, mycircle/.style={
                     circle,
                     draw=black,
                     fill=gray,
                     fill opacity = 0.3,
                     text opacity=1,
                     inner sep=0pt,
                     minimum size=20pt,
                     font=\small},
                  myarrow/.style={-Stealth},
                  node distance=0.6cm and 1.2cm]
    

                \draw[ultra thick, red, ->] (8.5,6.5) arc (120:80:3) ;
    
                \node at (13,9) [mycircle,draw] (j1) {Bus $j_1$};
                \node at (9,5) [red, dashed, circle, draw] (i) {Bus $i$};
                \node at (12.5,4.5) [mycircle, dashed, blue, draw] (j3) {Bus $j_3$};
                \draw[dashed, blue, -] (i) -- (j1);
                \draw[-] (j1) -- (j3);
                \draw[dashed, blue, -] (i) -- (j3);
    
                \node at (7.5,4.5) [mycircle,draw] (g1) {$g_1$};
                \node at (10.5,5.5) [mycircle, dashed, blue, draw] (g2) {$g_2$};
                \draw[ -] (i) -- (g1);
                \draw[ blue, -] (i) -- (g2);
    
                \node at (13.8,7.8) [mycircle,draw] (g3) {$g_3$};
                \draw[ -] (j1) -- (g3);
                \node at (7.8,8.8) (T) {{$I_i^\omega = \{i,\  g_2,\ (i, j_1),\ (i, j_3),\ j_3\}$. }};
                \node at (9,7.3) (T) {{fire spread direction}};
            \end{tikzpicture}
            \caption{A random set $I_i^\omega$ illustration. A fault occurs at bus $i$ (energized), which causes an endogenous wildfire to spread and affect two transmission lines $(i,j_1)$ and $(i,j_2)$, bus $j_3$, and generator $g_2$, which are marked by blue dashed lines.}
            \label{fig:endogenousfire}
        \end{figure}
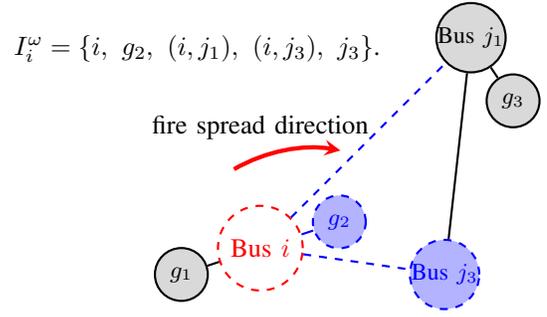

        \noindent We categorize wildfire risk into two classes: endogenous and exogenous wildfires. Exogenous wildfires result from external factors beyond the control of grid operators, while endogenous wildfires arise due to component faults. Fig.~\ref{fig:endogenousfire} is an illustration of a group of components affected by an endogenous fire.

        We conceptualize wildfire uncertainty as a stochastic disruption characterized by random factors, including its timing, location, and spread. To represent a single instance of such a disruption, we employ a portfolio of random variables indexed by $\omega \in \Omega$. Within this framework, we model i) the disruption's timing as $\tau^\omega$; ii) binary parameters $u^\omega_c$ and $v^\omega_c$ to account for both endogenous and exogenous wildfires; iii) the random set $I_c^\omega$ which describes the potentially affected components attributed to the endogenous wildfire event caused by component $c$. In the presence of multiple wildfire disruptions, each subsequent disruption is interdependent with the preceding one and shares a consistent structure.
       
    \subsection{Multistage Stochastic Mixed-Integer Programming Model} \label{subsec:formulation}
        \noindent When a fault occurs at an energized component, the associated endogenous wildfire may incur high damage costs. The power system operator can de-energize power equipment to reduce the endogenous risk of igniting a wildfire. However, suppose multiple components are de-energized to prevent endogenous fires. In that case, the power supply capacity may be greatly reduced, resulting in the inability to meet crucial load demand and causing serious secondary disasters. To improve the resilience of the power system and mitigate wildfire damage, operators should de-energize some potentially dangerous electrical equipment under high-risk conditions to ensure enough power supply while significantly reducing risk. 
    
        In model~\eqref{prob:stage1}, a de-energized generator will have zero generation capacity, and a de-energized line will be considered open. If a bus is shut off, all generators and lines connected to it will also be de-energized. Our model uses binary decision variables $z_{ct}$ to represent whether a component $c \in \cC$ is de-energized at time period $t$. We assume that a de-energized component will remain off until the end of the time horizon due to safety considerations \cite{rhodes2020balancing} unless it is restored. Model~\eqref{prob:stage1} obtains a nominal plan that should be implemented until a disruption occurs or the time horizon ends, whichever occurs first. If a fire disruption is observed at period $\tau^{\dot{\omega}}$ where $\dot{\omega} \in \Omega|_{\omega_0}$ and $\omega_0$ \tcr{represents the index of a deterministic realization for the $0$-th stage random variables}, given the current shut-off state $z_{\tau^{\dot{\omega}}-1}$, the model enters the disruptive stage and incurs the disruptive-stage value function $f^{\dot{\omega}}$, where we assume that the components at revealed ignition locations will be shut off.     
        \begin{subequations}
            \label{prob:stage1}
            \begin{align}
                & {Z^*(\beta) = \min \ \sum_{{\dot{\omega}} \in \Omega|_{\omega_0}} p^{\dot{\omega}} \left[ \sum_{t = 1}^{\tau^{\dot{\omega}}-1} \sum_{d \in \cD} w_d s_{dt} + f^{\dot{\omega}}(z_{\tau^{\dot{\omega}}-1}) \right] }\notag \\
                & \st \quad \forall t \in \cT: \notag \\
                & P^L_{ijt} \leq -b_{ij} \left(\theta_{it} -\theta_{jt} + \bar{\theta} (1-z_{ijt})\right) \quad\quad\ \ \forall (i,j) \in \cL \label{eqn:pfconsl1}\\
                & P^L_{ijt} \geq -b_{ij} \left(\theta_{it} -\theta_{jt} + \underbar{\theta} (1-z_{ijt})\right) \ \ \quad\quad\forall (i,j) \in \cL \label{eqn:pfconsr1}\\
                & -W_{ij} z_{ijt} \leq P^L_{ijt} \leq W_{ij} z_{ijt} \qquad\qquad\qquad \forall (i,j) \in \cL \label{eqn:thermallimit1}\\
                & {\sum_{g \in \cG_i} P^G_{gt} + \sum_{l \in \cL_i} P^L_{lt} = \sum_{d \in \cD_i} D_{dt} (1 - s_{dt}) \qquad\quad \forall i \in \cB} \label{eqn:flowbalance1}\\
                & \underbar{P}^G_g z_{gt} \leq P^G_{gt} \leq \bar{P}^G_g z_{gt} \quad\quad\quad\quad\quad\quad\quad\qquad\ \ \forall g \in \cG \label{eqn:genlimit1}\\
                & z_{it} \geq x_{dt} \quad\quad\quad\quad\quad\quad\quad\quad\quad\quad\quad\quad \ \ \forall i \in \cB, d \in \cD_i \label{eqn:loadlogic1}\\
                & z_{it} \geq z_{gt} \quad\quad\quad\quad\quad\quad\quad\quad\quad\quad\quad\quad\ \ \forall i \in \cB, g \in \cG_i \label{eqn:genlogic1}\\
                & z_{it} \geq z_{lt} \qquad\quad\quad\quad\quad\quad\quad\quad\quad\quad\quad\quad \forall i \in \cB, l \in \cL_i \label{eqn:linelogic1}\\
                & \tcr{r_{c,t-1} \geq r_{ct} \quad\qquad\quad\quad\quad\quad\quad\quad\quad\ \ \forall c \in \cC, \text{ if } t \geq 2} \label{eqn:restorationtime1} \\
                & \tcr{r_{c,t-1} - r_{ct} \geq z_{ct} - z_{c,t-1} \qquad\qquad\ \ \forall c \in \cC, \text{ if } t \geq 2} \label{eqn:componenttime1} \\
                & \sum_{t \in \cT} (s_{dt} - s_{d't}) \leq \beta \cdot T \qquad\qquad\qquad\qquad\ \forall d, d' \in \mathcal{D} \label{eqn:fairness}\\
                & z_{ct}, r_{ct} \in \{0,1\}\qquad\qquad\qquad\qquad\qquad\quad\quad\quad \forall c \in \cC.\label{eqn:binaryrestriction1}
            \end{align}
        \end{subequations}
        \tcr{For a given realization $\dot\omega$, the first-stage problem corresponds to the periods preceding the occurrence of the first disruption. Thus, the objective function considers the load-shedding cost from the initial period up to the period immediately preceding the disruption, and the expected cost $f^{\dot{\omega}}$ after the disruption $\dot{\omega}$ is observed.} Constraints~\eqref{eqn:pfconsl1}-\eqref{eqn:genlimit1} correspond to the DC power flow model and constraint~\eqref{eqn:loadlogic1} models the component interactions, equivalent to constraint (7)-(9) in \cite{rhodes2020balancing}. Combining constraint~\eqref{eqn:restorationtime1} with the binary restriction of $r_{ct}$ implies that the restoration operation $(r_{c1}, \dots, r_{cT})$ follows the pattern $(1, \dots, 1, 0, \dots, 0)$, with the first occurrence of $0$ indicating the restoration time, and a component $c \in \cC$ can only be restored once. We impose this constraint that each component can only be restored once. This restriction reflects the reality of short-term operations, where both starting up and shutting down a component \tcr{requires} time. It underscores the trade-off between shutting off and restoring a component.
        Constraints~\eqref{eqn:componenttime1} describe the temporal logic of components' status: once a component is shut off, it stays off unless restored. 
        Constraint~\eqref{eqn:fairness} ensures fairness by controlling the discrepancy between the maximum and minimum cumulative load-shedding percentages across loads. To modulate the stringency of fairness, we introduce a parameter $\beta$. A smaller value of $\beta$ enforces more rigorous fairness, ensuring that the difference in cumulative load-shedding percentages between any two loads does not exceed $\beta$. Importantly, this constraint maintains the mixed-integer programming (MIP) formulation, making it amenable to MIP solvers. 
        
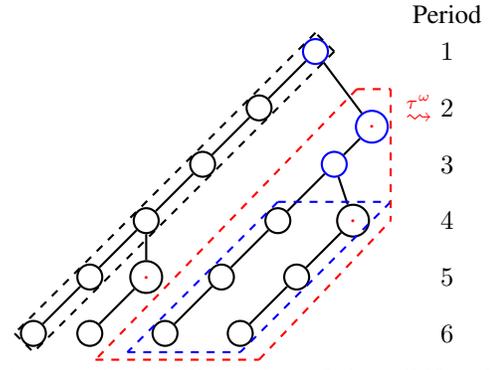
\begin{figure}
            \centering
            \begin{tikzpicture}[scale = .5, square/.style={regular polygon,regular polygon sides=4},node distance={20mm}, thick, main/.style = {draw, circle}, mycircle/.style={
                     circle,
                     draw=black,
                     fill=gray,
                     fill opacity = 0.3,
                     text opacity=1,
                     inner sep=0pt,
                     minimum size=20pt,
                     font=\small},
                  myarrow/.style={-Stealth},
                  node distance=0.6cm and 1.2cm]
    
                \node at (23.5,21) (Period) {Period};
                \node at (23.5,20) (Period1) {$1$};
                \node at (23.5,18.5) (Period2) {$2$};
                \node at (23.5,17) (Period3) {$3$};
                \node at (23.5,15.5) (Period4) {$4$};
                \node at (23.5,14) (Period5) {$5$};
                \node at (23.5,12.5) (Period5) {$6$};
                \node at (20,20) [circle, draw, blue] (R) {};
                \node at (18.5,18.5) [circle,draw] (N1) {};
                \node at (17,17) [circle, draw] (N2) {};
                \node at (15.5,15.5) [circle,draw] (N3) {};                        
                \node at (14,14) [circle,draw] (N4) {};
                \node at (12.5,12.5) [circle,draw] (N5) {};
                \node at (21.5,18) [circle,draw, blue] (11) {\textcolor{red}{.}};
                \node at (21,15.5) [circle,draw] (12) {\textcolor{red}{.}};
                \node at (20.5,17) [circle, draw, blue] (21) {};
                \node at (19.5,14) [circle, draw] (22) {};
                \node at (19,15.5) [circle,draw] (31) {};
                \node at (18,12.5) [circle,draw] (33) {};
                \node at (17.5,14) [circle,draw] (41) {};
                \node at (15.5,14) [circle,draw] (51) {\textcolor{red}{.}};
                \node at (16,12.5) [circle,draw] (52) {};
                \node at (14,12.5) [circle,draw] (61) {};
                
                \draw[-] (R) -- (11);
                \draw[-] (R) -- (N1);
                \draw[-] (N1) -- (N2);
                \draw[-] (N2) -- (N3);
                \draw[-] (N3) -- (N4);
                \draw[-] (N4) -- (N5);
                \draw[-] (11) -- (21);
                \draw[-] (21) -- (12);
                \draw[-] (12) -- (22);
                \draw[-] (21) -- (31);
                \draw[-] (31) -- (41);
                \draw[-] (22) -- (33);
                \draw[-] (N3) -- (51);
                \draw[-] (51) -- (61);
                \draw[-] (41) -- (52);
                
                \draw[-, dashed, black] (20.5,20) -- (12.5,12);
                \draw[-, dashed, black] (20,20.5) -- (12,12.5); 
                \draw[-, dashed, black] (20.5,20) -- (20,20.5);
                \draw[-, dashed, black] (12,12.5) -- (12.5,12);
    
                \draw[-, dashed, black, blue] (19,16) -- (22,16); 
                \draw[-, dashed, black, blue] (15,12) -- (18,12); 
                \draw[-, dashed, black, blue] (22,16) -- (18,12); 
                \draw[-, dashed, black, blue] (19,16) -- (15,12); 
    
                \draw[-, dashed, black, red] (21.1,19) -- (22,19); 
                \draw[-, dashed, black, red] (14.15,11.8) -- (18.5,11.8); 
                \draw[-, dashed, black, red] (22,19) -- (22,15.5); 
                \draw[-, dashed, black, red] (21.1,19) -- (14.15,11.8); 
                \draw[-, dashed, black, red] (18.5,11.8) -- (22,15.5); 
         
    
    
                \node at (22.75,18.5) [red] {{$\overset{\tau^\omega}{\rightsquigarrow}$}};
            \end{tikzpicture}
            \caption{A scenario-tree illustration of the wildfire disruption problem with $T = 6$.  The presence of a ``\textcolor{red}{.}'' symbol within the nodes denotes the occurrence of a disruption.
    }
            \label{fig:scenariotree}
        \end{figure}

        The scenario tree of the M-SMIP model is shown as Fig.~\ref{fig:scenariotree}, in which each node represents the decisions at a time period, and paths from the root to leaves represent scenarios. The black box of the left diagonal branch represents the nominal plan obtained if no disruption occurs, and each branch extending to the right represents a disruption scenario. Suppose we have a first-stage realization $\omega \in \Omega|_{\omega_0}$ with $\tau^\omega = 2$. The red box on the right corresponds to problem~\eqref{prob:staget} with $\tau^\omega = 2$ and some $\omega \in \Omega|_{\omega_0}$, of which the value function is $f^\omega$. For realization $\omega$, it has two scenario paths. One is a nominal scenario, and the other is a disruptive scenario, with the second disruption occurring at period $\tau^{\dot{\omega}} = 4$. The branches in the blue box are the scenarios with the same disruption in period $2$ and share the same historical information (the nodes in blue). 

        \tcr{For a given disruptions, we model the development of the post-disruption system in a similar manner as the nominal scenario as in the initial state. Essentially, we develop a new model that represent the post-disruption nominal plan. }The value function \(f^{\omega}\) for the $s$-th stage incorporates the shut-off state variables $\hat{z}^{\omega'}_{\tau^{\omega}-1}$, where $\omega'$ represents a realization of stage $s-1$ and $\omega \in \Omega|_{\omega'}$, and is characterized by the wildfire disruption realization $\xi^{\omega}$. It focuses on the operations after \(\tau^{\omega}\) and can be evaluated by the subproblem $f^\omega$:
        \begin{subequations}
            \label{prob:staget}
            \begin{align}
                & f^{\omega}(\hat{z}^{\omega'}_{\tau^{\omega}-1}) = \notag\\
                & \min \sum_{\dot{\omega}\in\Omega|_{\omega}} p^{\dot{\omega}} \left[\sum_{t = \tau^{\omega}}^{\tau^{\dot{\omega}}-1} \sum_{d \in \cD} w_d s^{\omega}_{dt} + \sum_{c \in \cC} c^r_c \nu^{\omega}_{c} + f^{\dot{\omega}}({z}^{\omega}_{\tau^{\dot{\omega}} - 1})\right] \notag\\
                & \st \quad \forall t \in \{\tau^{\omega},\dots, T\}: \notag \\
                & P^{L,{\omega}}_{ijt} \leq -b_{ij} \left(\theta^{\omega}_{it} -\theta^{\omega}_{jt} + \bar{\theta} (1-z^{\omega}_{ijt}) \right) \qquad \forall (i,j) \in \cL\label{eqn:pfconslt}\\
                & P^{L,{\omega}}_{ijt} \geq -b_{ij} \left(\theta^{\omega}_{it} -\theta^{\omega}_{jt} + \underbar{\theta} (1-z^{\omega}_{ijt})\right) \qquad \forall (i,j) \in \cL \label{eqn:pfconsrt}\\
                & -W_{ij} z^{\omega}_{ijt} \leq P^{L,{\omega}}_{ijt} \leq W_{ij} z^{\omega}_{ijt}\qquad\qquad\quad\ \ \forall (i,j) \in \cL \label{eqn:thermallimitt}\\
                & {\sum_{g \in \cG_i} P^{G,{\omega}}_{gt} + \sum_{l \in \cL_i} P^{L,{\omega}}_{lt} = \sum_{d \in \cD_i} D_{dt} (1 - s^{\omega}_{dt}) \quad\  \forall i \in \cB} \label{eqn:flowbalancet}\\
                & \underbar{P}^{G,{\omega}}_g z^{\omega}_{gt} \leq P^{G,{\omega}}_{gt} \leq \bar{P}^{G,{\omega}}_g z^{\omega}_{gt} \qquad\qquad\qquad\quad\ \ \forall g \in \cG \label{eqn:genlimitt}\\
                & z^{\omega}_{it} \geq x^{\omega}_{dt} \ \qquad\qquad\qquad\qquad\qquad\qquad\ \forall i \in \cB, d \in \cD_i \label{eqn:loadlogict}\\
                & z^{\omega}_{it} \geq z^{\omega}_{gt} \ \qquad\qquad\qquad\qquad\qquad\qquad\ \forall i \in \cB, g \in \cG_i \label{eqn:genlogict}\\
                & z^{\omega}_{it} \geq z^{\omega}_{lt} \qquad\qquad\qquad\qquad\qquad\qquad\quad \forall i \in \cB, l \in \cL_i \label{eqn:linelogict}\\
                & z^{\omega}_{ct-1} \geq z^{\omega}_{ct} \qquad\qquad\qquad\qquad\qquad\qquad\qquad\quad\forall c \in \cC \label{eqn:operationlogict}\\
                & z^{\omega}_{ct} \leq 1 - \nu_c^{\omega}\qquad\qquad\qquad\qquad\qquad\quad \qquad\quad\ \forall c \in \cC \label{eqn:fireFunction} \\
                & \nu_c^{\omega} \geq v_c^{\omega} \qquad\qquad\qquad\qquad\qquad\qquad\qquad\qquad \forall c \in \cC \label{eqn:exogenousFire} \\
                & \nu^{\omega}_k \geq u_c^{\omega} z^{\omega}_{c\tau^{\omega}-1} \qquad\qquad\qquad\qquad\quad\quad \forall c \in \cC, k \in I^{\omega}_{c} \label{eqn:ignitiont} \\
                & z^{\omega}_{c\tau^{\omega}-1} = \hat{z}^{\omega'}_{c \tau^{\omega} - 1}\qquad\qquad\qquad\qquad\qquad\qquad\ \forall c \in \cC \label{eqn:nonanticipativity}\\
                & z^{\omega}_{ct}, \nu_c^{\omega}, z^{\omega}_{ct} \in \{0,1\} \qquad\qquad\qquad\qquad\qquad\quad\ \forall c \in \cC. \label{eqn:binaryrestriction}
            \end{align}
        \end{subequations}
        \tcr{The objective function retains a similar meaning to the previous one, with the addition of a new term specifically introduced to account for the damage cost.} The wildfire damage cost for a component \(c \in \cC\), denoted by \(c_c^r\), consists of the replacement cost of the electric components and the financial loss to the nearby communities. In model~\eqref{prob:staget}, the energization status will stay the same in the remaining time horizon, as we assume that all wildfire damages reveal at period \(\tau^{\omega_t}\) and no recovery decisions take place afterward. Constraints~\eqref{eqn:pfconslt}-\eqref{eqn:genlimitt} model an equivalent form of the DC power flow constraints as in model~\eqref{prob:stage1}. Constraint~\eqref{eqn:loadlogict} indicates the functioning state of components, similar to their nominal-stage counterpart~\eqref{eqn:loadlogic1}. We create a local copy of the last-stage shut-off decisions, \(z^{\omega}_{c\tau^{\omega}-1}\), via the duplicating constraint~\eqref{eqn:nonanticipativity}, which is known as the nonanticipativity constraint. With \(z^{\omega}_{ct}\) indicating whether component \(c\) has been shut off, constraint~\eqref{eqn:operationlogict} makes sure that the shut-off components remain shut off during the second stage. Constraint~\eqref{eqn:fireFunction} states that damaged components no longer function, where we model the exogenous fire damage by constraint~\eqref{eqn:exogenousFire} and the endogenous fire damage by constraint~\eqref{eqn:ignitiont}. Notice that an endogenous fire started at component \(c\) requires both a fault occurrence \(u_c^{\omega_t} = 1\) and the component not being shut off \(z_c^{\omega_t} = 1\) and spreads to components \(k \in I_c^{\omega_t}\).
    
        We posit that at every stage $s$, a nominal realization $\omega^0_s$ exists, implying the absence of new disruptions. This realization lacks offspring realizations, i.e., $\Omega|_{\omega^0} = \emptyset$, and its associated value function equals $0$.

\section{Solution Methods}\label{sec:alg}
    \noindent A standard multistage stochastic program with time-period uncertainty can be solved using decomposition algorithms like the SDDP algorithm~\cite{pereira1991}. Our problem is different in setting that a `stage' includes all decisions between two disruptions, and its length is a random variable.  This section introduces an algorithm for solving problem~\eqref{prob:stage1}.
    \subsection{Decomposition Algorithm}
        \noindent The value function $f^{\omega}$ is non-convex nature and lack of an analytical expression. A widely employed approach for approximating value functions is the utilization of cutting planes. We adopt a collection of cuts to establish lower approximations for each value function $f^{\omega}$. Feasibility cuts are not required since the relatively complete recourse property always applies. This can be demonstrated by setting all variables in problem~\eqref{prob:staget} to zero. Specifically, for the problem~\eqref{prob:stage1} with a fixed $\beta$, we define its cutting-plane lower approximation as follows:
        \begin{align}\label{prob:relaxationstage1}
             \underline{Z}_\ell^* = \min &\ \sum_{\dot{\omega} \in \Omega|_{\omega_0}} p^{\dot{\omega}} \left[ \sum_{t = 1}^{\tau^{\dot{\omega}}-1} \sum_{d \in \cD} w_d s_{dt} + V^{\dot{\omega}} \right] \tag{$M_\ell$} \\
             \st & \ \mbox{Constrictions}~\eqref{eqn:pfconsl1}-\eqref{eqn:binaryrestriction1} \qquad\qquad\quad \forall t \in \cT \notag \\
             & \quad V^{\dot{\omega}} \geq (\lambda^{\dot{\omega}, k})^\top (z_{\tau^{\dot{\omega}} - 1} - \hat{z}_{\tau^{\dot{\omega}} - 1}^{k} ) + \notag \\
             & \qquad\qquad v^{\dot{\omega}, k}, \quad \forall \dot{\omega} \in \Omega|_{\omega_0}, k = 1,\dots, \ell-1, \notag
        \end{align}
        where the value function of the subsequent realization $\omega$, denoted as $f^\omega$, is approximated using $\ell-1$ cuts. Similarly, for the value function $f^{\omega}$, we define its cutting-plane lower approximation $\underbar{f}^{\omega}$ as follows:
        \begin{align}\label{prob:relaxationstaget}
             \underbar{f}^{\omega}_\ell & (\hat{z}^{\omega',\ell}_{\tau^{\omega}-1}) = \tag{$S_\ell^\omega$}\\
             \min \ & \sum_{\dot{\omega}\in\Omega|_{\omega}} p^{\dot{\omega}} \left[\sum_{t = \tau^{\omega}}^{\tau^{\dot{\omega}}-1} \sum_{d \in \cD} w_d s^{\omega}_{dt} + \sum_{c \in \cC} c^r_c \nu^{\omega}_{c} + V^{\dot{\omega}} \right] \notag\\
             \st \ & \mbox{Constrictions}~\eqref{eqn:pfconslt}-~\eqref{eqn:ignitiont} , ~\eqref{eqn:binaryrestriction} \quad \forall t \in \{\tau^{\omega},\dots, T\}  \notag\\
            & V^{\dot{\omega}} \geq (\lambda^{\dot{\omega}, k})^\top (z_{\tau^{{\dot{\omega}}} - 1}^{\omega} - \hat{z}_{\tau^{\dot{\omega}} - 1}^{\dot{\omega}, k} ) +\notag\\
            & \qquad\qquad\qquad\ v^{\dot{\omega}, k}, \qquad \forall \dot{\omega} \in \Omega|_{\omega}, k = 1,\dots, \ell-1 \notag\\
            & z^{\omega}_{\tau^{\omega}-1} = \hat{z}^{\omega',\ell}_{\tau^{\omega} - 1}, \tag{Nonanticipativity $:\lambda$}
        \end{align}
        During the $\ell$-th iteration of Algorithm~\ref{alg:decomposition}, an extra cut is generated to characterize $f^\omega$ by solving a relaxation problem of $\underbar{f}^{\omega}$ at the point $\hat{z}_{\tau^\omega - 1}^{\omega, \ell}$. This procedure yields the cut, which is subsequently incorporated into~\eqref{prob:relaxationstaget}, featuring the cut intercept $v^{\omega, \ell}$ and the cut slope $\lambda^{\omega,\ell}$. We refer to problem~\eqref{prob:relaxationstaget} that has additionally incorporated cuts as \textit{augmented} $(S_\ell^\omega)$. Without loss of generality, we can assign the disruption time for the nominal realization $\omega^0$ as $T+1$. As a result, the corresponding value function for $\omega^0$ equals $0$.

        \begin{algorithm2e}
            \small
            \SetAlgoLined
            Initialization cut iteration number $\ell = 1$, lower bound $LB = 0$, upper bound $UB = \infty$ and $\epsilon \geq 0$\; 
            \While{$\frac{UB-LB}{UB} \geq \epsilon$}{
                \tcc{Forward Steps}
                Solve problem~\eqref{prob:relaxationstage1} the first-stage shut-off solution $\hat{z}^\ell$, optimal value $\underline{Z}_\ell^*$, and the approximations of value function $\hat{V}^{\omega,\ell}$ for each $\omega \in \Omega|_{\omega_0}$ \;
                Update $LB = \underline{Z}_\ell^*$\;
                Let $\overline{Z} \leftarrow \underline{Z}_\ell^* - \sum_{\omega \in \Omega|_{\omega_0}} p^\omega \hat{V}^{\omega,\ell}$ and $p \leftarrow 1$\;
                \For{$\omega \in \Omega|_{\omega'}$}{
                    Solve problem~\eqref{prob:relaxationstaget} with $\hat{z}^{\omega',\ell}_{\tau^\omega - 1}$ to obtain the corresponding shut-off solution $\hat{z}^{\omega,\ell}$, optimal value $\underline{f}^{\omega}_\ell(\hat{z}^{\omega',\ell}_{\tau^\omega - 1})$, and the approximations of value function $\hat{V}^{\dot{\omega},\ell}$ for each $\dot{\omega} \in \Omega|_{\omega}$\;
                    $\overline{Z} \leftarrow \overline{Z} + p[\underline{f}^\omega_\ell(\hat{z}^{\omega',\ell}_{\tau^\omega - 1}) - \sum_{\dot{\omega} \in \Omega|_{\omega}} p^{\dot{\omega}} \hat{V}^{\dot{\omega},\ell}]$\;
                    \If{$\Omega|_{\omega} \not= \emptyset $}{
                        \For{$\dot{\omega} \in \Omega|_{\omega}$}{
                            $p \leftarrow p\cdot p^{\omega}$,
                            $\omega' \leftarrow \omega$,
                            $\omega \leftarrow \dot{\omega}$\;
                            Go to line $6$\;
                        }
                    }
                }
                \If{$\overline{Z} < UB$}{
                    Update $UB = \overline{Z}$ and shut-off solutions $\{\hat{z}^{\omega, \ell}\}_{\omega}$\;
                }
                \tcc{Backward Steps}
                \For{$\omega \in \Omega|_{\omega'}$}{
                    \If{$\Omega|_{\omega} \not= \emptyset $}{
                        \For{$\dot{\omega} \in \Omega|_{\omega}$}{
                            Go to line $19$\;
                        }
                    }
                    Solve a relaxed problem $({R}^\omega_\ell)$ of the \textit{augmented} $(S^\omega_\ell)$ to obtain the cut slope $\lambda^{\omega, \ell}$ and intercept $v^{\omega, \ell}$\; 
                    \textit{Augment} $(S_\ell^{\omega'})$ by adding the cut\;
                }
                Let $\ell = \ell + 1$\;
            }
            \textbf{Output:} The $\epsilon$-optimal value $UB$ and solutions $\{\hat{z}^{\omega, \ell}\}_{\omega}$.
            \caption{Decomposition algorithm based on Cutting-plane Method}\label{alg:decomposition}
        \end{algorithm2e}

    \subsection{Cut Families}
        \noindent This section covers a range of cut types and corresponding relaxed problem $(R^\omega_\ell)$ applicable in Algorithm~\ref{alg:decomposition}. 
        \subsubsection{Benders' Cut (BC)}
            The relaxation problem $({R}^\omega_\ell)$ solved in the backward step is the LP relaxation of the \textit{augmented} $({S}^\omega_\ell)$.
            Therefore, the cut coefficient $(v^{\omega, \ell}, \lambda^{\omega, \ell})$ is computed based on the optimal value of the LP relaxation and an optimal dual solution. Specifically, the cut added to its preceding realization $\omega'$ takes the following coefficient: $v^{\omega, \ell} = \underline{f}_\ell^{LP, \omega}(\hat{z}^{\omega',\ell}_{\tau^{\omega - 1}})$ and $\lambda^{\omega, \ell}$ is a basic optimal dual solution $\overline{\lambda}$ corresponding to the nonanticipativity constraint $z^{\omega}_{\tau^{\omega}-1} = \hat{z}^{\omega', \ell}_{\tau^{\omega}-1}$.

        \subsubsection{Lagrangian Cut (LC)}
            These cutting planes are derived by solving the Lagrangian dual relaxation problem~\eqref{opt:lagrangiandual} to acquire the optimal solution $\lambda^*$ for the cut slope $\lambda^{\omega, \ell}$ and the optimal value $\mathcal{R}_{\ell}^\omega(\hat{z}^{\omega',\ell}, \lambda^*; \mathcal{Z})$ for the intercept $v^{\omega, \ell}$, i.e., its relaxed problem $(R_\ell^\omega)$ is
            \begin{align}\label{opt:lagrangiandual}
                 \max_{\lambda} \quad \mathcal{R}_{\ell}^\omega(\hat{z}^{\omega',\ell}, \lambda; \mathcal{Z}),
            \end{align}
            where the Lagrangian relaxation problem $\mathcal{R}_{\ell}^\omega(\hat{z}^{\omega',\ell}, \lambda; \mathcal{Z})$ is derived by relaxing the nonanticipativity constraint and constraining the copy variables to the domain $\mathcal{Z}$ within the \textit{augmented} $(S_{\ell}^\omega)$ framework:
            \vspace{-.0cm}
            \begin{subequations}
                \label{prob:lagrangianrelaxation}
                \begin{align*}
                     \mathcal{R}_{\ell}^\omega & (\hat{z}^{\omega',\ell}, \lambda; \mathcal{Z}) = \\
                     \min \ & \sum_{\dot{\omega}\in\Omega|_{\omega}} p^{\dot{\omega}} \left[\sum_{t = \tau^{\omega}}^{\tau^{\dot{\omega}}-1} \sum_{d \in \cD} w_d s^{\omega}_{dt} +  \sum_{c \in \cC} c^r_c \nu^{\omega}_{c} + V^{\dot{\omega}} \right]  + \\
                     & \qquad\qquad\qquad\qquad\qquad\lambda^\top(\hat{z}^{\omega',\ell}_{\tau^{\omega} - 1} - z^{\omega}_{\tau^{\omega}-1})\\
                     \st \ & \mbox{Constrictions}~\eqref{eqn:pfconslt}-~\eqref{eqn:ignitiont} , ~\eqref{eqn:binaryrestriction} \quad \forall t \in \{\tau^{\omega},\dots, T\}  \notag\\
                     & V^{\dot{\omega}} \geq (\lambda^{\dot{\omega}, k})^\top (z_{\tau^{{\dot{\omega}}} - 1}^{\omega} - \hat{z}_{\tau^{\dot{\omega}} - 1}^{\dot{\omega}, k} ) + \\
                     & \qquad\qquad\qquad\qquad v^{\dot{\omega}, k}, \qquad \forall \dot{\omega} \in \Omega|_{\omega}, k = 1,\dots, \ell \\
                    & z_{c\tau^\omega - 1} \in \mathcal{Z} \qquad\qquad\qquad\qquad\qquad\qquad\qquad \forall c \in \cC.
                \end{align*}
            \end{subequations}

            \begin{remark}
                The standard Lagrangian cut method, as proposed by Zou et al.~\cite{zou2016nested}, involves selecting $\mathcal{Z} = [0,1]$. Alternatively, opting for $\mathcal{Z} = \{0,1\}$ yields improved cut performance.  
            \end{remark}
    
        \subsubsection{Strengthened Benders' Cut (SBC)}
            Given a fixed cut slope $\overline{\lambda}$, a valid cut can be derived, featuring a cut intercept of $\mathcal{R}_{\ell}^\omega (\hat{z}^{\omega',\ell}, \overline{\lambda}; \mathcal{Z})$. Subsequently, Benders' cut can be fortified through the aforementioned approach. The process involves initiating the dual solution $\overline{\lambda}$ corresponding to the nonanticipativity constraint for the LP relaxation of the \textit{augmented} $(S_\ell^\omega)$. The acquired dual solution $\overline{\lambda}$ then serves as the cut slope. It is followed by the resolution of a mixed-integer program to attain the strengthened cut intercept of $\mathcal{R}_{\ell}^\omega (\hat{z}^{\omega',\ell}, \overline{\lambda}; \mathcal{Z})$. In Fig.~\ref{fig:cut}, strengthened Benders' cut (in orange) is parallel with Benders' cut (in blue) while exhibiting improved performance.

        \subsubsection{Square-Minimization Cut (SMC)}
\begin{figure}[h]
                \begin{center}
                    \begin{tikzpicture}
                        \begin{axis}[
                                axis lines = middle,
                                xlabel = {$z$},
                                ylabel = {$f^\omega$},
                                xmin = -.08, xmax = 1.1,
                                ymin = 1.2, ymax = 3.3,
                                major grid style = {lightgray},
                                minor grid style = {lightgray!25},
                                width = .45\textwidth,
                                height = 0.35\textwidth]
                                \addplot[name path = Cut1, 
                                        domain = -0.1:1.1,
                                        dashed,
                                        very thick,
                                        blue,]{ 1.7 - 1 * x};
                                \addplot[name path = Cut3, 
                                        domain = -0.1:1.,
                                        dashed,
                                        very thick,
                                        black,]{ 3.0 - 5 * x};
                                \addplot[name path = Cut4, 
                                        domain = -0.1:1.2,
                                        dashed,
                                        very thick,
                                        orange,]{ 2.3 - 1 * x};
                                \addplot[name path = Cut5, 
                                        domain = -0.1:1.2,
                                        dashed,
                                        very thick,
                                        orange,]{ 2.1 - 1 * x};
                                \addplot[name path = Cut2, 
                                        domain = -0.1:1.1,
                                        dashed,
                                        very thick,
                                        red,]{ 3 - 1.7 * x};
                                    \addplot[name path = Cut6, 
                                        domain = -0.1:1.1,
                                        dashed,
                                        very thick,
                                        red,]{ 3 - 2.5 * x};
                                \addplot[mark=*,only marks, black] coordinates {(0, 1.2) };
                                \addplot[mark=*,only marks, teal] coordinates {(0,3) (1,1.3)};
                                \node[black] at (-.05,1.3) {$\hat{z}$};
                                \draw [dashed] (axis cs:{1.0,0}) -- (axis cs:{1.0,3.5});
                                \legend{BC, LC, SBC, , SMC}
                            \end{axis}
                        \end{tikzpicture}
                    \end{center}
                \caption{SMC and LC are tight at \(\hat{z} = 0\), whereas BC and SBC are not, as their values at $\hat{z}$ are smaller than the function value. The upper and lower SBC (SMC) are obtained by taking $\mathcal{Z} = \{0,1\}$ and $[0,1]$, respectively.}\label{fig:cut}
            \end{figure}
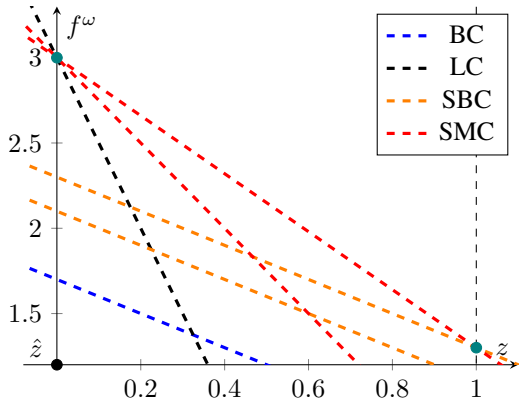

            Lagrangian cuts, represented by the black line in Fig.~\ref{fig:cut}, may be steep and fail to provide a good lower approximation at other solutions. To address this limitation, Square-Minimization Cut (SMC)~\cite{yang2023multi}, depicted by the red lines in Fig.~\ref{fig:cut}, offers a rotated approach. Instead of solving problem~\eqref{opt:lagrangiandual} to obtain cut coefficients, we use an alternative cut-generation problem~\eqref{prob:SMCgeneration} as the relaxed problem $(R_\ell^\omega)$:
            \vspace{-.cm}
            \begin{subequations}\label{prob:SMCgeneration}
                \begin{align}
                   \min_{\lambda}\quad & \lambda^\top \lambda \\
                   \mbox{s.t.}\quad & \mathcal{R}_{\ell}^\omega(\hat{z}^{\omega',\ell}, \lambda; \mathcal{Z}) \geq (1 - \delta)\underline{f}_\ell^{\omega}(\hat{z}^{\omega',\ell}_{\tau^{\omega - 1}}), \label{eqn:smc_cons}
                \end{align}
            \end{subequations}
            and we let \(\lambda^{\omega,\ell}\) equal to its optimal solution $\lambda^*$ and $v^{\omega,\ell} = \mathcal{R}_{\ell}^\omega(\hat{z}^{\omega',\ell}, \lambda^*; \mathcal{Z})$. Fig.~\ref{fig:cut} shows that the difference between a steep cut and a flat cut lies in their angle with the horizontal plane. We prefer a flat cut, corresponding to a smaller \(\lambda^\top \lambda\), as \(\lambda\) represents the linear cut's coefficients. We set up constraint~\eqref{eqn:smc_cons} to force the cut value to be within a \(\delta\) neighborhood of \(f^\omega\) at \(\hat{z}\), which can be considered an ``anchor point." As the function $\mathcal{R}_{\ell}^\omega$ is a concave function of \(\lambda\) given $\hat{z}^{\omega',\ell}$, constraint~\eqref{eqn:smc_cons} characterizes a convex set. We can use convex programming solution methods to solve problem~\eqref{prob:SMCgeneration}.

    \subsection{Convergence}
        \noindent We present the detailed steps of the decomposition algorithm in Algorithm~\ref{alg:decomposition}, which iteratively updates the bounds. In the end, we terminate Algorithm~\ref{alg:decomposition} once the relative gap is within a predefined tolerance threshold $\epsilon \geq 0$. We can show that Algorithm~\ref{alg:decomposition} converges to the optimal value in the finite step.
        \begin{theorem}[Convergence]
            When $\epsilon = 0$, Algorithm~\ref{alg:decomposition} terminates in a finite number of iterations and outputs an optimal solution to problem~\eqref{prob:stage1} after finitely many iterations, if the backward steps generate Lagrangian cuts or square-minimization cuts with $\delta = 0$.
        \end{theorem}

\section{Numerical Results}\label{sec:numerical}
    \noindent In this section, we outline the configuration of our numerical experiments employed to assess the efficacy of our decomposition method using various cut families. We examine the fairness of the M-SMIP formulation~\eqref{prob:stage1} across different fairness levels and restoration strategies, thereby demonstrating the advantages of M-SMIP.

    \subsection{Experiment Setup}
\begin{figure}
            \centering
            \setkeys{Gin}{width=.85\linewidth}
            \includegraphics{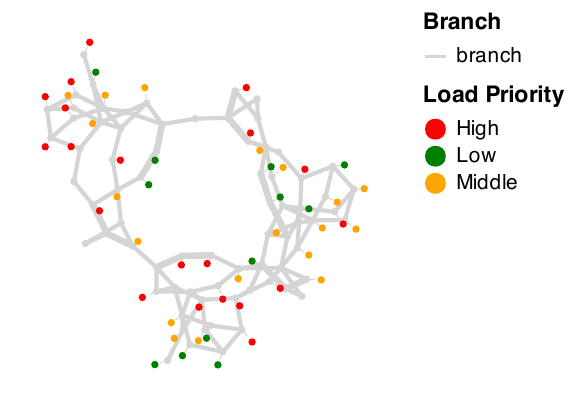}
                \caption{Illustration of the loads in RTS-GMLC system.}
                \label{fig:configuration}
        \end{figure}

        \noindent We utilize the RTS-GMLC 73-bus case~\cite{RTS-GLMC} for our study and conduct experiments over a $24$-hour horizon, with each time period encompassing three hours $(T = 8)$. The economic consequences of wildfires and sudden power interruptions are contingent upon both fire intensity and the load magnitude within the affected area. \tcr{To capture the economic impacts of wildfires and power outages, we evaluate the significance of each electrical component and its effect on the surrounding area using the same cost parameters as in Ref.~\cite{yang2023multi}:
        \begin{enumerate} 
            \item Load priorities \(w_d\) range from $50$ to $1000$; 
            \item The damage costs \(r_c\) of wind turbines, thermal and nuclear power plants are $50$, $1000$, and $2500$, respectively; 
            \item The damage cost \(r_c\) of each bus is $50$; 
            \item The damage cost \(r_c\) of a transmission line is $0.285 \ell$, where $\ell$ is the length of the transmission line. 
        \end{enumerate}
        Fig.~\ref{fig:configuration} color codes the load priority levels based on their load-shedding and damage costs, respectively. }
        
        All optimization models were implemented using JuMP~\cite{JuMP} in Julia v1.9~\cite{Julia} and solved by Gurobi v10.0.0~\cite{gurobi} on a computer with a $10$-core M1 Pro CPU and $32$GB memory. The network plots are generated using PowerPlot.jl, which depends on PowerModels.jl package~\cite{PowerModels}. Scenario simulation is constructed by an agent-based model package, Agent.jl~\cite{Agents}. For Algorithm~\ref{alg:decomposition} and SMC, we set $\epsilon = 1\%$ and $\delta = 10^{-4}$. 

    \subsection{Scenario Generation}
        \noindent We employ the term ``scenario'' to describe a complete disruption path within the scenario tree, extending from the root node to one of its leaves. For instance, there are four distinct scenarios in Fig.~\ref{fig:scenariotree}. We used identical settings and parameters to the cellular automaton simulation model described in Ref.~\cite{yang2023multi} but allowed for the possibility of multiple disruptions. 
        
        A \textit{training set}, $\Xi$, consists of $500$ samples \tcr{ which has been proven to be sufficient in Ref.~\cite{yang2023multi}.} A \textit{testing set}, $\tilde{\Omega}$, consists of $5,000$ samples, each of which contains at most one disruption. \tcr{Each scenario is assigned an equal probability of occurrence. Vulnerable components, identified as risky due to their heightened susceptibility to wildfire attacks, are associated with a greater likelihood of damage in multiple scenarios. These risky components receive heightened attention within the sample sets, reflecting their increased probability of being affected.}


    \subsection{Cut Performance}
\begin{figure}
            \centering
            \setkeys{Gin}{width=.95\linewidth}
            \includegraphics{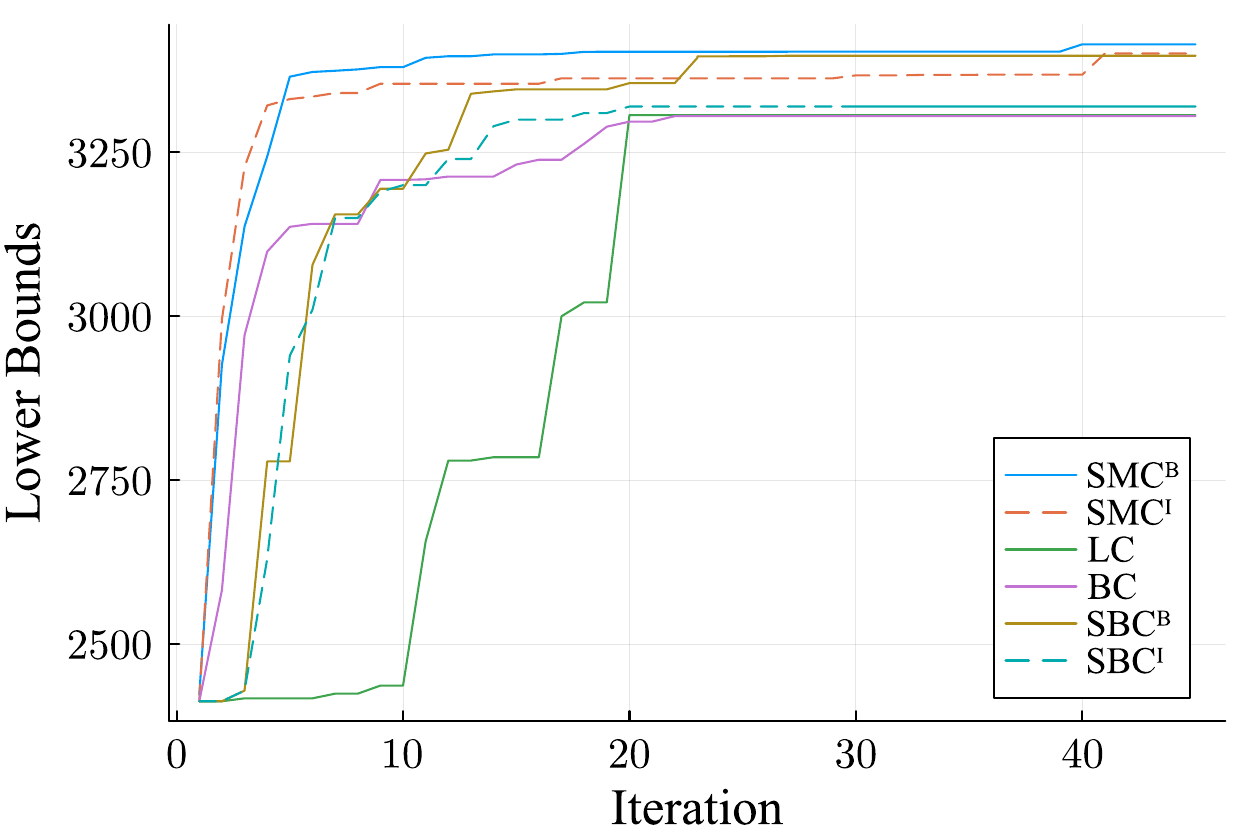}
                \caption{Lower Bounds vs. Iteration Count. Superscripts are employed to denote the choice of $\mathcal{Z}$, where $B$ represents the binary set, $\{0,1\}$, and $I$ denotes the interval, $[0,1]$.}
                \label{fig:convergenceiter}
        \end{figure}

        \noindent We applied Algorithm~\ref{alg:decomposition} using four cut families to the scenario set $\Xi$ with $\beta = 0.4$ and assessed their performance. Fig.~\ref{fig:convergenceiter} depicts the relationship between convergence and algorithmic iteration count across different cut families. 

        Our findings suggest that:
        i) Benders' cuts fail to achieve convergence due to their lack of tightness;
        ii) Square-minimization cuts converge with fewer iterations, although each iteration requires more computation time compared to Lagrangian cuts;
        iii) Cuts with $\mathcal{Z} = \{0,1\}$ exhibit superior performance in terms of lower-bound improvement compared to those with $[0,1]$.

        In TABLE~\ref{table:time}, a maximum runtime of $25,000$ seconds is specified. This table presents the best bounds, gap, time per iteration, and the total time required to reach these bounds for each cut type, signifying the minimum duration necessary. BC and LC cannot achieve convergence within the time limit. SBC can be considerably tighter due to their improved intercept achieved through MIP solving. Two versions of SMCs underscore the advantages of selecting $\mathcal{Z} = \{0,1\}$. 

        \begin{table}[h]
            \centering
            \caption{Assessing Algorithm~\ref{alg:decomposition} performance with different cut families.}
            \begin{tabular}{c|ccc|cc}\toprule[1pt]\hline
               \multirow{2}{*}{Cut Type} & \multicolumn{3}{c|}{Solution Quality}  & \multicolumn{2}{c}{Runtime (sec.)} \\              & LB                  & UB              & Gap          & Time/Iter.  & Total  \\ 
                \cmidrule{1-6} 
            BC        & $3305.3$             & $3457.8$        & $4.41\%$     & $214.2$     & $5377$                     \\
            LC        & $3306.9$             & $4116.8$        & $19.67\%$    & $193.4$     & $3803$                     \\
            SBC$^B$   & $3397.5$             & $3449.6$        & $1.51\%$     & $237.7$     & $8472$                     \\
            SBC$^I$   & $3323.1$             & $3433.6$        & $3.15\%$     & $210.8$     & $10500$                     \\
            SMC$^B$   & $3414.9$             & $3426.5$        & $0.34\%$     & $445.4$     & $20516$                    \\
            SMC$^I$   & $3368.6$             & $3419.3$        & $1.48\%$     & $391.9$     & $13414$                    \\\hline  
            \bottomrule[1pt]    
            \end{tabular}
            \label{table:time}
        \end{table}


    \subsection{Fairness Evaluation}\label{sec:fairness}
\begin{figure*}
            \centering
            \setkeys{Gin}{width=\linewidth}
            \begin{subfigure}{0.32\textwidth}
                \includegraphics{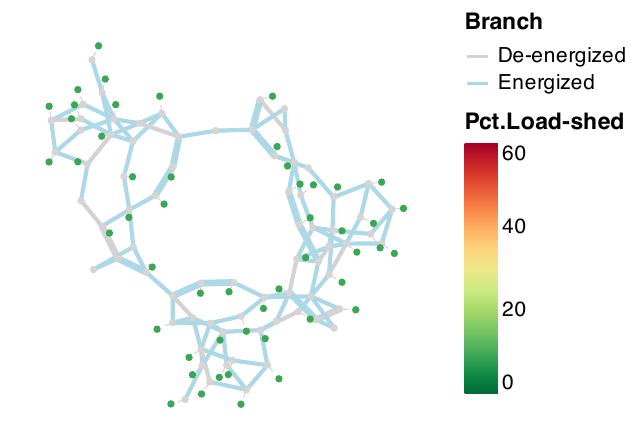}
                \caption{When $\beta = 0.0$, load-shedding is evenly distributed across numerous loads, with a maximum load-shedding percentage of $10\%$, resulting in a total load-shedding of $9.2\%$.}
                \label{fig:fairness1}
            \end{subfigure}
            \hfil
            \begin{subfigure}{0.32\textwidth}
                \includegraphics{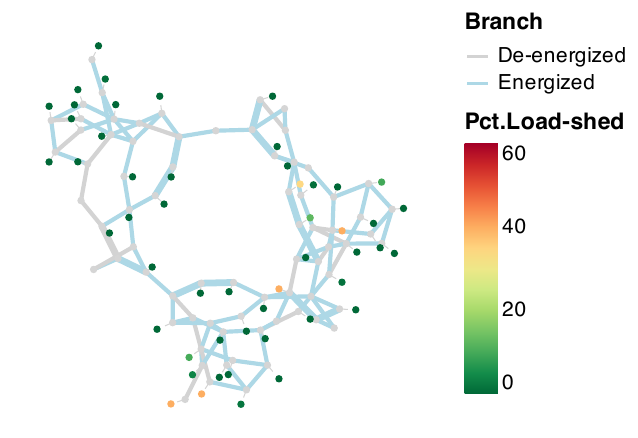}
                \caption{When $\beta = 0.4$, load-shedding is less evenly distributed, with a maximum load-shedding percentage of $40\%$, resulting in a total load-shedding of $4.5\%$. }
                \label{fig:fairness3}
            \end{subfigure}
            \hfil
            \begin{subfigure}{0.32\textwidth}
                \includegraphics{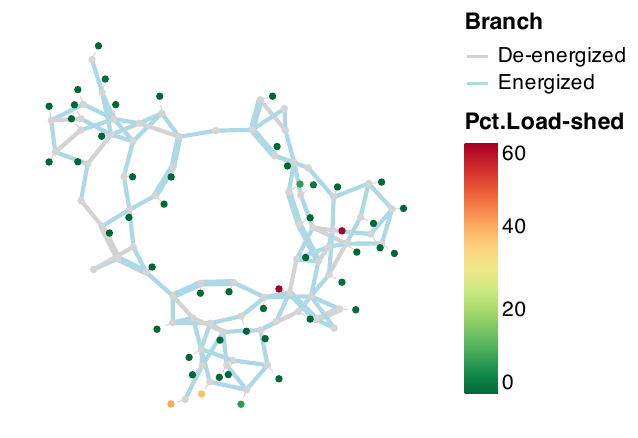}
                \caption{When $\beta = 0.6$, load-shedding is concentrated on a few lower-priority loads, with a maximum load-shedding percentage of $60\%$, resulting in a total load-shedding of $4.2\%$. }
                \label{fig:fairness6}
            \end{subfigure}
            \caption{Percentage of load-shedding for each load with different $\beta$ and de-energization operations illustration.}
            \label{fig:fairness}
        \end{figure*}

        \noindent We examine different values of $\beta$ to evaluate the trade-off between efficiency and fairness. 
        \tcr{A solution that neglects or minimally considers fairness may lead to a significant imbalance in load shedding among different regions, as illustrated in Fig.~\ref{fig:fairness6}. On the other hand, a solution that enforces strict fairness may be too conservative and incur a high cost, as shown in Fig.~\ref{fig:fairness1}. Our demonstration reveals that a balanced and efficient solution can be achieved with a moderate level of fairness, exemplified in Fig.~\ref{fig:fairness3}.}
        The divergence in total load-shedding or the emergence of more uniform load-shedding can be attributed to changes in power distribution. 
        In contrast to the approach taken in the work~\cite{Kody2022}, 
        our de-energization decisions are more consistent and allow you to select your desired level of fairness.

    \subsection{Comparative Analysis: Restoration vs.~No Restoration}
        \noindent In a prior study~\cite{yang2023multi}, we explored the benefits of stochastic programming, demonstrating its superiority by achieving a minimum improvement of at least $45\%$ compared to deterministic models. In this work, we shift our focus to restoration analysis. Given $\beta$, solving problem~\eqref{prob:stage1} with the scenario set $\Xi$ yields a nominal shutoff plan, denoted as $X^{*,\beta} = \{s^*, z^*, \theta^*, P^*\}.$ To gauge its effectiveness, we compare this plan with alternatives $X^{\prime,\beta} = \{s', z', \theta', P'\}$ generated without considering restoration by replacing the restoration constraints~\eqref{eqn:restorationtime1} and~\eqref{eqn:componenttime1} with the following \textit{component-time logic constraint} $z_{c,t-1} \geq z_{ct}$.
        We evaluate the expected total cost for each SAA solution (nominal plan) using the scenario set $\tilde{\Xi}$ for the out-of-sample test as follows: 
        \tcr{
            \begin{align*}
                & g(\hat{X}) = \min \ \sum_{\dot{\omega} \in \tilde{\Omega}} \frac{1}{|\tilde{\Omega}|} \left[ \sum_{t = 1}^{\tau^{\dot{\omega}}-1} \sum_{d \in \cD} w_d \hat{s}_{dt} + f^{\dot{\omega}}(\hat{z}_{\tau^{\dot{\omega}}-1}) \right]. \label{eqn:obj}
            \end{align*}
        Nominal restoration plans tend to de-energize more components than non-restoration plans, as they can restore certain components once they are no longer at risk. This is advantageous in wildfire scenarios as it reduces damage costs, but it may increase load-shedding costs in the nominal scenario due to more de-energized components.}
        
        We present the results in TABLE~\ref{table:restorationtest}, displaying the benefits of integrating restoration operations and shedding light on the performance of nominal plans across different fairness levels denoted by $\beta$. When \(\beta > 0\), these restoration operations yield more than 5\% cost reductions as they enhance the capacity for restoring components, allowing for their de-energization when necessary. Subsequently, these components can be restored after hazardous periods, resulting in a minor reduction in nominal load supply but a considerable decrease in total costs.

        As detailed in Section~\ref{sec:fairness}, varying values of $\beta$ can result in different nominal plans that demonstrate similar disruptive performance, confirming their relatively consistent de-energization operations. However, smaller $\beta$ values necessitate a considerably higher level of additional load-shedding.

        \begin{table}[ht]
            \centering
            \caption{Nominal load-shedding cost, disruptive load-shedding costs, and damage costs under different nominal plans obtained using different settings.}
            \label{table:restorationtest}
            \begin{tabular}{cc|c|cc|c}\toprule[1pt]\hline
                \multicolumn{2}{c|}{Setting}  & \multicolumn{1}{c|}{Nominal}                      & \multicolumn{2}{c|}{Disruptive}    & Total Cost     \\ 
                 Res.               & $\beta$ & Load shed                  & Load shed             & \multicolumn{1}{c|}{Damage}        & $g(\cdot)$        \\ \cmidrule{1-6} 
\multirow{4}*{$X^{*,\beta}$}   & $0.0$   & $23276.1$                  & $1789.0$              & $1389.0$                           & $18517.6$          \\
                                    & $0.2$   & $3974.7$                   & $1840.0$              & $1457.4$                           & $6554.2$          \\
                                    & $0.4$   & $2978.4$                   & $1971.7$              & $1515.6$                           & $5948.9$          \\
                                    & $0.6$   & $2873.5$                   & $1967.9$              & $1533.7$                           & $5886.4$          \\\cmidrule{1-6}
    \multirow{4}*{$X^{\prime,\beta}$}    & $0.0$   & $21710.1$                  & $2108.8$              & $1957.4$                           & $18370.4$         \\ 
                                    & $0.2$   & $3819.0$                   & $1947.8$              & $1863.7$                           & $6957.4$          \\ 
                                    & $0.4$   & $2933.7$                   & $2077.3$              & $1887.3$                           & $6409.6$          \\ 
                                    & $0.6$   & $2776.7$                   & $2045.8$              & $1937.2$                           & $6331.7$          \\\hline 
            \bottomrule[1pt]    
            \end{tabular}
        \end{table}

\section{Conclusion}\label{sec:conclusion}

    \noindent This work represents a multistage stochastic mixed-integer program for power system operations under the persistent wildfire threat. Our multistage stochastic program is adept at comprehensively modeling the intricate uncertainties stemming from spatially varying wildfire disruptions and their temporal evolution. Our model exhibits robustness and resilience in addressing fairness concerns and indicates the benefit of adding a restoration option. 

    To complement our modeling approach, we propose an efficient decomposition algorithm that capitalizes on binary state variables, enabling the generation of valid cuts and enhancing the resolution of extensive instances. Empirical assessments underscore the effectiveness of various cut families, highlighting the overall efficiency of our multistage model. In future work, we aim to develop theory-driven cutting planes for our decomposition algorithm, balance computational efficiency, and integrate unit commitment decisions with ramping constraints. \tcr{AC power flow equations can be incorporated in the formulation via additional constraints or approximations, and it is necessary to explore how they affect the solution's fairness and robustness in future work.}

\bibliographystyle{IEEEtran}
\bibliography{PSCC_Ref}

\end{document}